\documentclass[12pt]{article}
\usepackage{amsthm,amssymb}
\def\bbr{{\mathbb R}}
\def\esssup{\mathop{\rm ess\; sup}}
\def\essinf{\mathop{\rm ess\; inf}}
\theoremstyle{definition}
\newtheorem{exam}{Example}[section]
\newtheorem{dfn}[exam]{Definition}
\newtheorem*{rem}{Remark}
\theoremstyle{plain}
\newtheorem{thm}[exam]{Theorem}

\begin{document}
\title{
Bayesian Posteriors Without Bayes' Theorem}

\author{T. P. Hill and Marco Dall' Aglio}
\date{}
\maketitle
\begin{abstract}
The classical Bayesian posterior arises naturally as the
unique solution of several different optimization problems, without the
necessity of interpreting data as conditional probabilities and then
using Bayes' Theorem. For example, the classical Bayesian posterior is
the unique posterior that minimizes the loss of Shannon information in
combining the prior and the likelihood distributions. These results,
direct corollaries of recent results about conflations of probability
distributions, reinforce the use of Bayesian posteriors, and may help
partially reconcile some of the differences between classical and
Bayesian statistics.
\end{abstract}

\section{Introduction}\label{sec1}

In statistics, prior belief about the value of an unknown parameter,
$\theta\in\Theta\subseteq\bbr^n$
obtained from experiments or other methods, is often expressed as a Borel
probability distribution $P_0$ on 
$\Theta\subseteq\bbr^n$ called the {\it prior distribution}. New evidence
or information about the value of $\theta$, based on an independent experiment
or survey, is recorded as a {\it likelihood distribution} $L$.  Here and
throughout it will be assumed that the likelihood function has finite
positive total mass, and that $L$ has been normalized, so that in fact 
$L$ is also
a Borel probability distribution on $\Theta$. 
Given the prior distribution  $P_0$ and
the likelihood distribution $L$, a {\it posterior distribution}
$P_1=P_1(P_0,L)$  for $\theta$ incorporates
the new likelihood information about $\theta$ into the information from the 
prior,
thus updating the prior. The posterior distribution
$P_1$   is typically viewed
as the conditional distribution of $\theta$
 given the new likelihood information,
often expressed as a random variable $X$.

The first main goal of this note is to use recent results for conflations
of probability distributions \cite{3,4} to show that the Bayesian posterior
is the unique posterior that minimizes the loss of Shannon information in
combining the prior and likelihood distributions. The Bayesian posterior
is also the unique posterior that attains the minimax likelihood ratio
of the prior and likelihood distributions, and the unique posterior
that is a proportional consolidation of the prior and likelihood
distributions. Thus, the classical Bayesian posterior appears naturally
as the solution of several different optimization problems, without the
necessity of interpreting likelihood as a conditional probability and
then invoking  Bayes Theorem. These results reinforce the use of Bayesian
posteriors, and may help partially reconcile some of the differences
between classical statistics and Bayesian statistics.

The second main goal of this note, another direct corollary of recent
results for conflations of probability distributions 
\cite{4},  is to identify
the best posterior when the prior and likelihood distributions are not
weighted equally, such as in cases when the prior distribution is given
more weight than the likelihood distribution. This new weighted posterior,
the unique distribution that minimizes the loss of weighted Shannon
information, coincides with the classical Bayesian posterior if the
prior and likelihood are weighted equally, but in general is different.

\section{Combining Priors and Likelihoods into Posteriors}\label{sec2}

There are many different methods for combining several probability
distributions (e.g., see \cite{1,3}), and in particular, for combining the
prior distribution  $P_0$ and the likelihood distribution $L$ into a single
{\it posterior distribution} $P_1=P_1(P_0,L)$. 
For example, the prior and likelihoods could
simply be averaged, i.e.\ $P_1=\frac{P_0+L}2$, 
or the data underlying the prior and the
likelihood could be averaged, in which case the posterior $P_1$ would be
the distribution of $\frac{X_0+X_L}2$, 
where $X_0$ and $X_L$ are independent random variables with
distributions $P_0$ and $L$, respectively.

In Bayesian statistics, the likelihood function $L$ is usually interpreted
as $L(\theta)=\alpha P(X\mid \theta)$,
where $X$ is the independent experiment or random variable yielding new
information about $\theta$, and $\alpha$ is the normalizing constant for 
$L$ to have mass one
(cf.\ \cite{2}). The 
{\it Bayesian posterior distribution} $P_B$ is then calculated using
Bayes Theorem: for example, if $P_0$ and $L$  are discrete with probability mass
functions (p.m.f.'s) $p_0$ and  $p_L$  respectively, then 
$P_B$ is discrete with p.m.f.
$$
p_B(\theta)=\frac{p_0(\theta)p_L(\theta)}{\sum_{\hat\theta\in\Theta}
p_0(\widehat\theta)p_L(\widehat\theta)}\,;$$
and if $P_0$ and  $L$ are absolutely continuous with probability density 
functions
(p.d.f.'s) $f_0$ and  $f_L$ respectively, then $P_B$ is absolutely continuous 
with p.d.f.
$$
f_B(\theta)=\frac{f_0(\theta)f_L(\theta)}{\int_\Theta
f_0(\widehat\theta)f_L(\widehat\theta)d\widehat\theta}$$
(provided the denominators are positive and finite).

\section{Minimizing Loss of Shannon Information}\label{sec3}

When the goal is to consolidate information from a prior distribution
and a likelihood distribution into a (posterior) distribution, replacing
those two distributions by a single distribution will clearly result
in some loss of information, however that is defined. Recall that
the classical {\it Shannon information} (also called the 
{\it self-information} or
{\it surprisal}) associated with the event $A$ from 
a probability distribution $P$,
$S_P(A)$, is given by $S_P(A)=-\log_2P(A)$
 (so the smaller the value of $P(A)$, the greater the information
or surprise). The numerical value of the Shannon information of a given
probability is simply the number of binary bits of information reflected
in that probability.

\begin{exam}\label{ex3.1} %%%Example 3.1. 
\rm If $P$ is uniformly distributed on $(0,1)$ and 
$A=(0,0.25)\cup (0.5,0.75)$, then $S_P(A)=-\log_2(P(A))=-\log_2(0.5)=1$, 
so	if $X$ is
a random variable with distribution $P$, then exactly one binary bit of
information is obtained by observing that $X\in A$, in this case that the value
of the second binary digit of $X$ is $0$.
\end{exam}

\begin{dfn}\label{dfn3.2} %%%Definition 3.2. 
\rm The {\it combined Shannon Information} associated with the
event $A$ from the prior distribution $P_0$ and the likelihood distribution  
$L$
is 
$$S_{\{P_0,L\}} (A)=S_{P_0}(A)+S_L(A)=-\log_2P_0(A)L(A),$$
 and the {\it maximum loss between the Shannon Information of a posterior
distribution $P_1$ and the combined Shannon information} of the prior and
likelihood distributions $P_0$ and $L$, $M(P_1;P_0,L)$, is
$$
M(P_1;P_0,L)=\max_A \left\{S_{\{P_0,L\}} (A)-S_{P_1}(A)\right\}=
\max_A\left\{\log_2 \frac{P_1(A)}{P_0(A)L(A)}\right\}.
$$
\end{dfn}
Note that the definition of combined Shannon information implicitly
assumes independence of the prior and likelihood distributions.  Note also
that no information is obtained by observing an event that is certain
to occur, so for instance $S_{[P_0,L]}(\Theta)=S_{P_1}(\Theta)=0$.
This implies that $M(P_1;P_0,L)$ is never negative.

\begin{dfn}\label{dfn3.3} %%%Definition 3.3. 
A prior distribution $P_0$ and a likelihood distribution $L$ are
{\it compatible} if $P_0$ and $L$ are both discrete with p.m.f's $p_0$ and
$p_L$ satisfying $\sum_{\theta\in\Theta}p_0(\theta)p_L(\theta)>0$,
or are both absolutely continuous with 
p.d.f.'s $f_0$ and $f_L$ satisfying $0<\int_\Theta f_0(\theta)f_L(\theta)
d\theta<\infty$.
\end{dfn}

\begin{exam}\label{ex3.4} %%%Example 3.4. 
Every two geometric distributions are compatible, every two
normal distributions are compatible, and every exponential distribution
is compatible with every normal distribution. Distributions with disjoint
support, discrete or continuous, are not compatible.
\end{exam}

\begin{rem} In practice, compatibility is not problematic. Any two
distributions may be easily transformed into two new distributions,
arbitrarily close to the original distributions, so that the two new
distributions are compatible, for instance by convolving each with
a  $U(-\epsilon,\epsilon)$ distribution.
\end{rem}

\begin{thm}\label{thm3.5} %%%Theorem 3.5. 
Let $P_0$ and $L$ be discrete compatible prior and likelihood
distributions. Then the Bayesian posterior $P_B$ is the unique posterior
distribution that minimizes the maximum loss of Shannon information from
the prior and likelihood distributions, i.e., that minimizes 
$M(P_1;P_0,L)$ among all
posterior distributions $P_1$. Moreover,
$$M(P_1;P_0,L)\ge \log_2\left[\left(\sum_{\theta\in\Theta}
p_0(\theta)p_L(\theta)\right)^{-1}\right]\mbox{ for all posterior
distributions } P_1,$$
 and equality is uniquely attained by the Bayesian posterior 
$P_1=P_B$.
\end{thm}

The conclusion of Theorem~\ref{thm3.5} follows immediately as a special case of
\cite[Corollary~4.4]{3}; analogous conclusions for the case of compatible
absolutely continuous distributions follow from \cite[Theorem~4.5]{3}. For
the benefit of the reader, a sketch of the proof of Theorem~\ref{thm3.5} similar
to that in \cite{4} is included.

\noindent{\it Sketch of proof.} 
First observe that for an event $A$, the difference between
the combined Shannon information obtained from a prior distribution
$P_0$ and
a likelihood distribution $L$, and the Shannon information obtained from
the posterior $P_1$, is
$$S_{\{P_0,L\}}(A)-S_{P_1}(A)=S_{P_0}(A)+S_L(A)-S_{P_1}(A)=\log_2
\frac{P_1(A)}{P_0(A)L(A)}\,.$$
Since $\log_2(x)$ is strictly increasing, the maximum (loss) thus occurs 
for an
event $A$ where 
$\frac{P_1(A)}{P_0(A)L(A)}$ is maximized.  

Next note that the largest loss of Shannon
information occurs for small sets $A$, since for disjoint sets $A$ and $B$,
$$\frac{P_1(A\cup B)}{P_0(A\cup B)L(A\cup B)}\le
\frac{P_1(A)+P_1(B)}{P_0(A)L(A)+P_0(B)L(B)}\le\max
\left\{\frac{P_1(A)}{P_0(A)L(A)}\,, \frac{P_1(B)}{P_0(B)L(B)}\right\},$$
where the inequalities follow from the inequalities  
$(a+b)(c+d)\ge ac+bd$ and  $\frac{a+b}{c+d}\le\max\left\{
\frac ac\,,\frac bd\right\}$ for positive
numbers $a,b,c,d$. Thus the problem reduces to finding the probability
mass function  $p$  that makes the maximum, over all real values 
$\theta$, of the ratio $\frac{p(\theta)}{p_0(\theta)p_L(\theta)}$
as small as possible.  But the minimum over all nonnegative
$q_1,\dots, q_n$ with $q_1+\cdots + q_n=1$
of the maximum of  $\frac{q_1}{r_1},\dots,\frac{q_n}{r_n}$
occurs when $\frac{q_1}{r_1}=\cdots=\frac{q_n}{r_n}$ 
(if they are not equal, reducing the
numerator of the largest ratio, and increasing that of the smallest,
will make the maximum smaller).  Thus the $p$ that makes the maximum of
$\frac{p(\theta)}{p_0(\theta)p_L(\theta)}$ as
small as possible is when $p(\theta)=cp_0(\theta)p_L(\theta)$, 
where $c$ is chosen to make $p$ a probability mass
function, i.e., to make $p(\theta)$ sum to 1. 
But this is exactly the definition of
the Bayesian posterior $P_B$ in the discrete case.~\hfill$\Box$

\section{Minimax Likelihood Ratios}\label{sec4}

In classical hypotheses testing, a standard technique to decide from
which of several known distributions given data actually came is to
maximize the likelihood ratios, that is, the ratios of the p.m.f.'s
or p.d.f.'s. Analogously, when the objective is to decide how best to
consolidate a prior distribution $P_0$ and a likelihood distribution $L$ into a
single (posterior) distribution 
$P_1=P_1(P_0,L)$, one natural criterion is to choose $P_1$ so
as to make the ratios of the likelihood of observing 
$\theta$ under $P_1$ as close as
possible to the likelihood of observing $\theta$ 
under {\it both} the prior distribution $P_0$ 
and the likelihood distribution  $L$. This motivates the following notion
of minimax likelihood ratio posterior.

\begin{dfn}\label{dfn4.1} %%%Definition 4.1. 
A discrete probability distribution $P^*$ (with p.m.f.\ $p^*$) is the
{\it minimax likelihood ratio (MLR) posterior} of  a discrete prior distribution
$P_0$ with p.m.f.\ $p_0$ and a discrete likelihood distribution  
$L$  with p.m.f.\ $p_L$ if
$$
\min_{\mbox{\scriptsize p.m.f.'s}\; p}\left\{\max_{\theta\in\Theta}
\frac{p(\theta)}{p_0(\theta)p_L(\theta)}-\min_{\theta\in\Theta}
\frac{p(\theta)}{p_0(\theta)p_L(\theta)}\right\}$$
is attained by $p=p^*$ (where 
$0/0:=1$).  
\end{dfn}

Similarly, an a.c.\ distribution $P^*$ with p.d.f.\ $f^*$ 
is the MLR posterior of an a.c.\ prior distribution $P_0$ with p.d.f.\ 
$f_0$ and an
a.c.\ likelihood distribution $L$ with p.d.f.\ $f_L$  if
$$
\min_{\mbox{\scriptsize p.m.f.'s}\; f}\left\{\esssup_{\theta\in\Theta}
\frac{f(\theta)}{f_0(\theta)f_L(\theta)}-\essinf_{\theta\in\Theta}
\frac{f(\theta)}{f_0(\theta)f_L(\theta)}\right\}$$
is attained by $f^*$.

The min-max terms in Definition~\ref{dfn4.1} are similar to the min-max criterion
for loss of Shannon Information (Theorem~\ref{thm3.5}), whereas the others are
dual max-min criteria. Just as the Bayesian posterior minimizes the loss
of Shannon information, the Bayesian posterior is also the MLR posterior
of the prior and likelihood distributions.

\begin{thm}\label{thm4.2} %%%Theorem 4.2. 
Let $P_0$ and $L$ be compatible discrete or compatible absolutely
continuous prior and likelihood distributions, respectively. Then the
unique MLR posterior for $P_0$ and $L$ is the Bayesian posterior distribution 
$P_B$.
\end{thm}

\noindent{\it Proof.} Immediate from \cite[Theorem~5.2]{3}.~\hfill$\Box$

\section{Proportional Posteriors}\label{sec5}

A criterion similar to likelihood ratios is to require that the posterior
distribution $P_1$  reflect the relative likelihoods of identical individual
outcomes under both $P_0$ and $L$. For example, if the probability that the
prior and the (independent) likelihood are both 
$\theta_a$ is twice that of the
probability both are $\theta_b$, then $P_1(\theta_a)$ 
should also be twice as large as $P_1(\theta_b)$.

\begin{dfn}\label{dfn5.1} %%%Definition 5.1. 
A discrete (posterior) probability distribution $P^*$ with
p.m.f.\ $p^*$ is a {\it proportional posterior of a discrete prior distribution
$P_0$ with p.m.f.\ $p_0$ and a compatible  discrete likelihood distribution	
$L$
with p.m.f.}\ $p_L$ if
$$\frac{p^*(\theta_a)}{p^*(\theta_b)}=\frac{p_0(\theta_a)p_L(\theta_a)}
{p_0(\theta_b)p_L(\theta_b)}\quad\mbox{for all } \theta_a,\theta_b\in
\Theta.$$
\end{dfn}

Similarly, a posterior a.c.\ distribution
$P^*$  with p.d.f.\ $f^*$ is a {\it proportional
posterior of an a.c.\ prior distribution $P_0$  with p.d.f.\
$f_0$ and a compatible
likelihood distribution $L$ with p.d.f.}\ $f_L$ if
$$\frac{f^*(\theta_a)}{f^*(\theta_b)}=\frac{f_0(\theta_a)f_L(\theta_a)}
{f_0(\theta_b)f_L(\theta_b)}\quad\mbox{ for (Lebesgue) almost all }
\theta_a,\theta_b\in \Theta.$$

\begin{thm}\label{thm5.2} %%%Theorem 5.2. 
Let $P_0$ and $L$ be compatible discrete or compatible absolutely
continuous prior and likelihood distributions, respectively. Then the
Bayesian posterior distribution $P_B$ is a proportional consolidation  for
$P_0$ and $L$.
\end{thm}

\noindent{\it Proof.} Immediate from \cite[Theorem 5.5]{3}.~\hfill$\Box$

\section{Optimal Posteriors for Weighted Prior and Likelihood Distributions}
\label{sec6}

\begin{dfn}\label{dfn6.1} %%%Definition 6.1.  
Given a prior distribution $P_0$ with weight $w_0>0$ and a likelihood
distribution  $L$ with weight $w_L>0$, the 
{\it combined weighted Shannon information}
associated with the event $A$, $S_{(P_0,w_0;L,w_L)}(A)$, is 
$$S_{(P_0,w_0;L,w_L)}(A)=\frac{w_0}{\max\{w_0,w_L\}}\,
S_{P_0}(A)+\frac{w_L}{\max\{w_0,w_L\}}\, S_L(A).$$
\end{dfn}

This definition ensures that only the {\it relative weights} are important,
so for instance if $w_0=w_L$, the combined weighted Shannon information of the
prior and likelihood always coincides with the (unweighted) combined
Shannon information of the prior and likelihood.  Note again that no
information is attained by observing any event that is certain to occur,
no matter what the distributions and weights, since 
$S_{P_0}(\Theta)=S_L(\Theta)=0$. The next theorem,
a special case of \cite[(8)]{4}, identifies the posterior distribution that
minimizes the loss of weighted Shannon information in the case the prior
and likelihood distributions are compatible discrete distributions;
the case for compatible absolutely continuous distributions is analogous.

\begin{thm}\label{thm6.2} %%%Theorem 6.2. 
Let $P_0$ and $L$ be compatible discrete prior and likelihood
distributions with p.m.f.'s $p_0$ and $p_L$ and weights
$w_0>0$ and $w_L>0$, respectively. Then the
unique posterior distribution that minimizes the maximum loss of Shannon
information from the weighted prior and likelihood distributions, i.e.,
that minimizes, among all posterior distributions $P_1$,
$$\max_A\left\{S_{(P_0,w_0;L,w_L)}(A)-S_{P_1}(A)\right\},$$
is the posterior distribution $P^w_1$ with p.m.f.
$$p^w_1(\theta)=\frac{(p_0(\theta))^{\frac{w_0}{\max[w_0,w_L]}}
(p_L(\theta))^{\frac{w_L}{\max[w_0,w_L]}}}
{\sum_{\hat\theta\in\Theta}(p_0(\widehat\theta))^{\frac{w_0}{\max[w_0,w_L]}}
(p_L(\widehat\theta))^{\frac{w_L}{\max[w_0,w_L]}}}\,.$$
\end{thm}

\begin{rem}  If both the prior and likelihood distributions are normally
distributed, the Bayesian posterior is also a best linear unbiased
estimator (BLUE) and a maximum likelihood estimator (MLE); e.g.\ see 
\cite{3}.
\end{rem}

\end{document}